\documentstyle[12pt]{article}
\title{The First Zagreb Index, the Forgotten Topological Index, the Inverse Degree and Some Hamiltonian Properties of Graphs}
\author {Rao Li\\
              Dept. of Computer Science, Engineering and Mathematics\\
              University of South Carolina Aiken\\
              Aiken, SC 29801 \\ 
              USA \\}
\date{Sept. 15, 2024}

\begin{document}
\maketitle
\begin{abstract}
Let $G = (V, E)$ be a graph. The first Zagreb index and the forgotten topological index of a graph $G$ are defined respectively as
$\sum_{u \in V} d^2(u)$ and $\sum_{u \in V} d^3(u)$, where $d(u)$ is  the degree of vertex $u$ in $G$. If the minimum degree of $G$ is at least one,
the inverse degree of $G$ is defined as $\sum_{u \in V} \frac{1}{d(u)}$. In this paper, we, for a graph with minimum degree at least one,
  present an upper bound for the first Zagreb index of the graph
and lower bounds for the forgotten topological index and the inverse degree of the graph. 
We also present sufficient conditions involving the first Zagreb index, the forgotten topological index, or the inverse degree for some Hamiltonian properties of a graph.
\end{abstract}  
$$Mathematics \,\, Subject \,\, Classification: 05C45, \,\, 05C09$$
$$Keywords:   The \,\, first \,\, Zagreb \,\, index, \,\,  the \,\, forgotten \,\, topological \,\, index,$$  
 $$the \,\, inverse \,\,  degree, \,\, Hamiltonian \,\, graph, \,\, traceable \,\, graph $$ 
\newpage 

 \noindent {\bf 1 \, Introduction } \\
 
\noindent We consider only finite undirected graphs without loops or multiple edges.
Notation and terminology not defined here follow those in \cite{Bondy}.
Let $G = (V(G), E(G))$ be a graph with $n$ vertices and $e$ edges, the degree of a vertex $v$ is denoted by $d_G(v)$. We use $\delta$ and $\Delta$ to denote the minimum degree and maximum degree of $G$, respectively.  
A set of vertices in a graph $G$ is independent if the vertices in the set are pairwise nonadjacent. A maximum independent set in a graph $G$ is an independent set of largest possible size. 
 The independence number, denoted $\beta(G)$, of a graph $G$ is the cardinality of a maximum independent set in $G$. 
 For disjoint vertex subsets $X$ and $Y$ of $V(G)$, we use 
 $E(X, Y)$ to denote the set of all the edges in $E(G)$ such that one end vertex of each edge is in $X$ and another end vertex of the edge is in $Y$. Namely, 
$E(X, Y) := \{\, e : e = xy \in E, x \in X, y \in Y \,\}$. 
A cycle $C$  in a graph $G$ is called a Hamiltonian cycle of $G$ if $C$ contains all the vertices of $G$. 
 A graph $G$ is called Hamiltonian if $G$ has a Hamiltonian cycle.
A path $P$  in a graph $G$ is called a Hamiltonian path of $G$ if $P$ contains all the vertices of $G$. 
A graph $G$ is called traceable if $G$ has a Hamiltonian path.\\

The first Zagreb index and the forgotten topological index  of a graph were introduced by Gutman and Trinajsti\'{c} \cite{GT} and Furtula and Gutman
\cite{FG}, respectively. For a graph $G$, its first Zagreb index,
denoted $Z_1(G)$, and its forgotten topological index, denoted $F(G)$, are defined as $\sum_{u \in V(G)} d_G^{2}(u)$ and
$\sum_{u \in V(G)} d_G^{3}(u)$, respectively. If $\delta(G) \geq 1$, the inverse degree, denoted $Inv(G)$, of $G$ is defined as $\sum_{u \in V} \frac{1}{d(u)}$.
  Using one inequality in \cite{D}, we in this paper present an upper bound for the first Zagreb index
and lower bounds for the forgotten topological index and the inverse degree of a graph $G$ with $\delta(G) \geq 1$.  
We also present sufficient conditions involving the first Zagreb index, the forgotten topological index, or the inverse degree for some Hamiltonian properties of a graph. The main results are as follows. \\

\noindent {\bf Theorem $1$.} Let $G$ be a graph with $n$ vertices, $e$ edges, and $\delta \geq 1$. Then \\

\noindent $[1]$ 
$$Z_1(G) \leq (n - \beta) \Delta^2 + \frac{e^2}{2 \beta} + \frac{\beta \Delta^3}{2 \delta}$$
with equality if and only if $G$ is a regular balanced bipartite graph. \\

\noindent $[2]$ 
$$F(G)  \geq (n - \beta) \delta^3 + \frac{\delta (2 \beta^2 \delta^2 - e^2)}{\beta}$$
with equality if and only if $G$ is a regular balanced bipartite graph. \\
 
\noindent $[3]$
$$F(G)  \geq (n - \beta) \delta^3 + \frac{\delta}{\beta} \left( 2 \beta \left(\beta \delta^2 + \frac{e^2}{n - \beta}\right) - e^2 - 2 \beta (n - \beta) \Delta^2\right)$$ \\
with equality if and only if $G$ is a regular balanced bipartite graph. \\

\noindent $[4]$ 
$$Inv(G) \geq \frac{n - \beta}{\Delta} + \frac{(2 \beta^2 \delta^2 - e^2)}{\beta \Delta^3}$$
with equality if and only if $G$ is a regular balanced bipartite graph. \\
 
\noindent $[5]$ 
$$Inv(G) \geq \frac{n - \beta}{\Delta} + \frac{1}{\beta \Delta^3} \left( 2 \beta \left(\beta \delta^2 + \frac{e^2}{n - \beta}\right) - e^2 - 2 \beta (n - \beta) \Delta^2\right)$$ 
with equality if and only if $G$ is a regular balanced bipartite graph. \\

\noindent {\bf Theorem $2$.} Let $G$ be a $k$-connected ($k \geq 2$) graph with $n \geq 3$ vertices and $e$ edges.  \\
 
\noindent [$1$] If 
 $$Z_1(G) \geq  (n - k - 1)\Delta^2 + \frac{e^2}{2(k + 1)} + \frac{(k + 1) \Delta^3}{2 \delta},$$
 then $G$ is Hamiltonian. \\ 
 
\noindent [$2$] If 
$$F(G) \leq (n - k - 1) \delta^3 + \frac{\delta (2 (k + 1)^2 \delta^2 - e^2)}{k + 1},$$ 
 then $G$ is Hamiltonian. \\ 
 
 \noindent [$3$] If 
 $$F(G) \leq (n - k - 1) \delta^3 + $$ $$\frac{\delta}{k + 1} \left( 2 (k + 1) \left((k + 1) \delta^2 + \frac{e^2}{n - k - 1}\right) 
 - e^2 - 2 (k + 1) (n - k - 1) \Delta^2\right),$$
 then $G$ is Hamiltonian. \\ 
 
\noindent [$4$] If
$$Inv(G) 
\leq \frac{n - k - 1}{\Delta} + \frac{(2 (k + 1)^2 \delta^2 - e^2)}{(k + 1) \Delta^3},$$ 
 then $G$ is Hamiltonian. \\ 
 
 \noindent [$5$] If 
 $$Inv(G) \leq \frac{n - k - 1}{\Delta} + $$
 $$\frac{1}{(k + 1) \Delta^3} \left( 2 (k + 1) \left((k + 1) \delta^2 + 
 \frac{e^2}{n - k - 1}\right) - e^2 - 2 (k + 1) (n - k - 1) \Delta^2\right),$$
 then $G$ is Hamiltonian. \\ 
 
 \noindent {\bf Theorem $3$.} Let $G$ be a $k$-connected ($k \geq 1$) graph with $n \geq 9$ vertices and $e$ edges.  \\
 
\noindent [$1$] If 
 $$Z_1(G) \geq  (n - k - 2)\Delta^2 + \frac{e^2}{2(k + 2)} + \frac{(k + 2) \Delta^3}{2 \delta},$$
 then $G$ is traceable. \\ 
 
\noindent [$2$] If 
$$F(G) \leq (n - k - 2) \delta^3 + \frac{\delta (2 (k + 2)^2 \delta^2 - e^2)}{k + 2},$$
 then $G$ is traceable. \\ 
 
 \noindent [$3$] If 
 $$F(G) \leq (n - k - 2) \delta^3 + $$ $$\frac{\delta}{k + 2} \left( 2 (k + 2) \left((k + 2) \delta^2 + \frac{e^2}{n - k - 2}\right) 
 - e^2 - 2 (k + 2) (n - k - 2) \Delta^2\right),$$
 then $G$ is traceable. \\ 
 
 \noindent [$4$] If 
 $$Inv(G) \leq \frac{n - k - 2}{\Delta} + \frac{(2 (k + 2)^2 \delta^2 - e^2)}{(k + 2) \Delta^3},$$
 then $G$ is traceable. \\ 
 
 \noindent [$5$] If 
 $$Inv(G) \leq \frac{n - k - 2}{\Delta} + $$
 $$\frac{1}{(k + 2) \Delta^3} \left( 2 (k + 2) \left((k + 2) \delta^2 + 
 \frac{e^2}{n - k - 2}\right) - e^2 - 2 (k + 2) (n - k - 2) \Delta^2\right),$$
 then $G$ is traceable. \\ 

\noindent {\bf 2 \, Lemmas } \\

We will use the following results as our lemmas. Lemma $1$ is Corollary $2.11$ on Page $8$ in \cite{D}. \\

\noindent {\bf Lemma $1$} \cite{D}. If $a_k$ and $b_k$ ($k = 1, 2, \cdots, s$) are positive real numbers, then
$$\frac{1}{2} \left(\sum_{i = 1}^s \frac{a_i^3}{b_i}\sum_{i = 1}^s \frac{b_i^3}{a_i} - \left(\sum_{i = 1}^s a_i b_i \right)^2 \right) \geq \sum_{i = 1}^s a_i^2 \sum_{i = 1}^s b_i^2 - \left(\sum_{i = 1}^s a_i b_i \right)^2 \geq 0.$$

The next two are from \cite{CE}. \\

\noindent {\bf Lemma $2$} \cite{CE}. Let $G$ be a $k$-connected graph of order $n \geq 3$. If $\beta \leq k$, then $G$ is Hamiltonian. \\

\noindent {\bf Lemma $3$} \cite{CE}. Let $G$ be a $k$-connected graph of order n. If $\beta \leq k + 1$, then $G$ is traceable. \\

Lemma $4$ below is from \cite{M}.\\

\noindent {\bf Lemma $4$} \cite{M}. Let $G$ be a balanced bipartite graph of order $2n$ with bipartition ($A$, $B$). If
$d(x) + d(y) \geq n + 1$ for any $x \in A$ and any $y \in B$ with $xy \not \in E$, then $G$ is Hamiltonian. \\

\noindent {\bf 3 \, Proofs } \\
 
\noindent {\bf Proof of Theorem $1$.} Let $G$ be a graph with $n$ vertices, $e$ edges, and $\delta \geq 1$. Clearly, $\beta < n$. 
 Let $I := \{\, u_1, u_2, ..., u_{\beta} \,\}$ be a maximum independent set in $G$. Then 
$$ \sum_{u \in I} d(u) = |E(I, V - I)| \leq \sum_{v \in V - I} d(v).$$ 
Since $\sum_{u \in I} d(u) +  \sum_{v \in V - I} d(v) = 2e$, we have that 
$$\sum_{u \in I} d(u) \leq e \leq \sum_{v \in V - I} d(v).$$
Applying Lemma $1$ with $s =  \beta$, 
$a_i = d(u_i)$ and $b_i = 1$ with $i = 1, 2, ... , \beta$, we have 
$$\frac{1}{2} \left(\sum_{i = 1}^{\beta} \frac{d^3(u_i)}{1}\sum_{i = 1}^s \frac{1^3}{d(u_i)} - \left(\sum_{i = 1}^{\beta} d(u_i)*1 \right)^2 \right) \geq \sum_{i = 1}^{\beta} d^2(u_i) \sum_{i = 1}^{\beta} 1^2 - \left(\sum_{i = 1}^{\beta} d(u_i)*1 \right)^2.$$
Thus
$$2 \beta \sum_{u \in I} d^2(u) \leq \sum_{u \in I} d^3(u) \sum_{u \in I} \frac{1}{d(u)} + \left(\sum_{u \in I} d(u) \right)^2
\leq \sum_{u \in I} d^3(u)  \sum_{u \in I} \frac{1}{d(u)} + e^2.$$
\noindent [$1$]. Notice that $Z_1(G) = \sum_{u \in I} d^2(u) + \sum_{v \in V - I} d^2(v)$. Thus 
$$2 \beta Z_1(G) = 2 \beta \sum_{u \in I} d^2(u) + 2 \beta \sum_{v \in V - I} d^2(v)$$
$$\leq \sum_{u \in I} d^3(u) \sum_{u \in I} \frac{1}{d(u)} + e^2
+ 2 \beta \sum_{v \in V - I} d^2(v) \leq \beta \Delta^3 \, \frac{\beta}{\delta} + e^2 + 2 \beta (n - \beta) \Delta^2.$$
Therefore 
$$Z_1(G) \leq (n - \beta) \Delta^2 + \frac{e^2}{2 \beta} + \frac{\beta \Delta^3}{2 \delta}.$$
If $$Z_1(G) = (n - \beta) \Delta^2 + \frac{e^2}{2 \beta} + \frac{\beta \Delta^3}{2 \delta},$$
then, from the above proofs, we have that  $\sum_{u \in I} d(u) = e$ which implies that $\sum_{v \in V - I} d(v) = e$ and thereby $G$ is a
bipartite graph with partition sets of $I$ and $V - I$. Furthermore, we have that $d(u) = \delta = \Delta$ for each $u \in I$ and  
$d(v) = \Delta$ for each $v \in I$. Thus $\delta |I| = |E(I, V - I)| = (n - |I|) \delta$. Therefore $|I| = |V - I| = \frac{n}{2}$.
Hence $G$ is a regular balanced bipartite graph. \\

\noindent If $G$ is a regular balanced bipartite graph, a simple computation can verify that  
$$Z_1(G) = (n - \beta) \Delta^2 + \frac{e^2}{2 \beta} + \frac{\beta \Delta^3}{2 \delta}.$$

This completes the proofs of [$1$]. \\

\noindent [$2$]. From 
$$2 \beta \sum_{u \in I} d^2(u) \leq \sum_{u \in I} d^3(u)  \sum_{u \in I} \frac{1}{d(u)} + e^2.$$
we have 
$$2 \beta \beta \delta^2 \leq \frac{\beta}{\delta} \sum_{u \in I} d^3(u) + e^2.$$
Thus $$\sum_{u \in I} d^3(u) \geq \frac{\delta (2 \beta^2 \delta^2 - e^2)}{\beta}.$$
Therefore $$F(G) = \sum_{w \in V} d^3(w) = \sum_{u \in I} d^3(u) + \sum_{v \in V - I} d^3(v)
\geq (n - \beta) \delta^3 + \frac{\delta (2 \beta^2 \delta^2 - e^2)}{\beta}.$$ 
If
$$F(G) = (n - \beta) \delta^3 + \frac{\delta (2 \beta^2 \delta^2 - e^2)}{\beta},$$
then, from the above proofs, we have that  $\sum_{u \in I} d(u) = e$ which implies that $\sum_{v \in V - I} d(v) = e$ and thereby $G$ is a
bipartite graph with partition sets of $I$ and $V - I$. Furthermore, we have that $d(u) = \delta$ for each $u \in I$ and  
$d(v) = \delta$ for each $v \in I$. Thus $\delta |I| = |E(I, V - I)| = (n - |I|) \delta$. Therefore $|I| = |V - I| = \frac{n}{2}$.
Hence $G$ is a regular balanced bipartite graph. \\

\noindent If $G$ is a regular balanced bipartite graph, a simple computation can verify that 
$$F(G) = (n - \beta) \delta^3 + \frac{\delta (2 \beta^2 \delta^2 - e^2)}{\beta}.$$ 

This completes the proofs of [$2$]. \\

\noindent [$3$]. By Cauchy-Schwarz inequality, we have
$$\sum_{v \in V - I} d^2(v) \sum_{v \in V - I} 1^2 \geq \left(\sum_{v \in V - I} d(v)\right)^2 \geq e^2.$$
Thus $$\sum_{v \in V - I} d^2(v) \geq \frac{e^2}{n - \beta}.$$
Since $$2 \beta \sum_{u \in I} d^2(u) \leq \sum_{u \in I} d^3(u)  \sum_{u \in I} \frac{1}{d(u)} + e^2,$$
we have 
$$2 \beta \left(\beta \delta^2 + \frac{e^2}{n - \beta}\right) \leq 2 \beta \sum_{u \in I} d^2(u) + 2 \beta \sum_{v \in V - I} d^2(v)$$
$$\leq \sum_{u \in I} d^3(u)  \sum_{u \in I} \frac{1}{d(u)} + e^2 + 2 \beta (n - \beta) \Delta^2 $$
$$\leq \frac{\beta}{\delta} \sum_{u \in I} d^3(u)  + e^2 + 2 \beta (n - \beta) \Delta^2.$$
Thus $$\sum_{u \in I} d^3(u) \geq \frac{\delta}{\beta} \left( 2 \beta \left(\beta \delta^2 + \frac{e^2}{n - \beta}\right) - e^2 - 2 \beta (n - \beta) \Delta^2\right).$$
Therefore 
$$F(G) = \sum_{w \in I} d^3(w) = \sum_{u \in I} d^3(u) + \sum_{v \in V - I} d^3(v)$$ 
$$\geq (n - \beta) \delta^3 + \frac{\delta}{\beta} \left( 2 \beta \left(\beta \delta^2 + \frac{e^2}{n - \beta}\right) - e^2 - 2 \beta (n - \beta) \Delta^2\right).$$
If 
$$F(G) = (n - \beta) \delta^3 + \frac{\delta}{\beta} \left( 2 \beta \left(\beta \delta^2 + \frac{e^2}{n - \beta}\right) - e^2 - 2 \beta (n - \beta) \Delta^2\right),$$
then, from the above proofs, we have that  $\sum_{u \in I} d(u) = e$ and $\sum_{v \in V - I} d(v) = e$ and thereby $G$ is a
bipartite graph with partition sets of $I$ and $V - I$. Furthermore, we have that $d(u) = \delta$ for each $u \in I$ and  
$d(v) = \delta = \Delta$ for each $v \in I$. Thus $\delta |I| = |E(I, V - I)| = (n - |I|) \delta$. Therefore $|I| = |V - I| = \frac{n}{2}$.
Hence $G$ is a regular balanced bipartite graph. \\

\noindent If $G$ is a regular balanced bipartite graph, a simple computation can verify that 
$$F(G) = (n - \beta) \delta^3 + \frac{\delta}{\beta} \left( 2 \beta \left(\beta \delta^2 + \frac{e^2}{n - \beta}\right) - e^2 - 2 \beta (n - \beta) \Delta^2\right).$$

This completes the proofs of [$3$]. \\

\noindent [$4$]. From 
$$2 \beta \sum_{u \in I} d^2(u) \leq \sum_{u \in I} d^3(u)  \sum_{u \in I} \frac{1}{d(u)} + e^2.$$
we have 
$$2 \beta \beta \delta^2 \leq \beta \Delta^3 \sum_{u \in I} \frac{1}{d(u)} + e^2.$$
Thus 
$$\sum_{u \in I} \frac{1}{d(u)} \geq \frac{(2 \beta^2 \delta^2 - e^2)}{\beta \Delta^3}.$$
Therefore 
$$Inv(G) = \sum_{w \in V} \frac{1}{d(w)} = \sum_{u \in I} \frac{1}{d(u)} + \sum_{v \in V - I} \frac{1}{d(v)}
\geq \frac{n - \beta}{\Delta} + \frac{(2 \beta^2 \delta^2 - e^2)}{\beta \Delta^3}.$$
If 
$$Inv(G) = \frac{n - \beta}{\Delta} + \frac{(2 \beta^2 \delta^2 - e^2)}{\beta \Delta^3},$$
then, from the above proofs, we have that  $\sum_{u \in I} d(u) = e$ which implies that $\sum_{v \in V - I} d(v) = e$ and thereby $G$ is a
bipartite graph with partition sets of $I$ and $V - I$. Furthermore, we have that $d(u) = \delta = \Delta$ for each $u \in I$ and  
$d(v) = \Delta$ for each $v \in I$. Thus $\Delta |I| = |E(I, V - I)| = (n - |I|) \Delta$. Therefore $|I| = |V - I| = \frac{n}{2}$.
Hence $G$ is a regular balanced bipartite graph. \\

\noindent If $G$ is a regular balanced bipartite graph, a simple computation can verify that  
$$Inv(G) = \frac{n - \beta}{\Delta} + \frac{(2 \beta^2 \delta^2 - e^2)}{\beta \Delta^3}.$$

This completes the proofs of [$4$]. \\

\noindent [$5$]. Recall that  $$\sum_{v \in V - I} d^2(v) \geq \frac{e^2}{n - \beta}.$$
Since $$2 \beta \sum_{u \in I} d^2(u) \leq \sum_{u \in I} d^3(u)  \sum_{u \in I} \frac{1}{d(u)} + e^2,$$
we have 
$$2 \beta \left(\beta \delta^2 + \frac{e^2}{n - \beta}\right) \leq 2 \beta \sum_{u \in I} d^2(u) + 2 \beta \sum_{v \in V - I} d^2(v)$$
$$\leq \sum_{u \in I} d^3(u)  \sum_{u \in I} \frac{1}{d(u)} + e^2 + 2 \beta (n - \beta) \Delta^2 $$
$$\leq \beta \Delta^3 \sum_{u \in I} \frac{1}{d(u)} + e^2 + 2 \beta (n - \beta) \Delta^2.$$
Thus $$\sum_{u \in I} \frac{1}{d(u)} \geq \frac{1}{\beta \Delta^3} \left( 2 \beta \left(\beta \delta^2 + \frac{e^2}{n - \beta}\right) - e^2 - 2 \beta (n - \beta) \Delta^2\right).$$
Therefore 
$$Inv(G) = \sum_{w \in V} \frac{1}{d(w)} = \sum_{u \in I} \frac{1}{d(u)} + \sum_{v \in V - I} \frac{1}{d(v)}$$
$$\geq \frac{n - \beta}{\Delta} + \frac{1}{\beta \Delta^3} \left( 2 \beta \left(\beta \delta^2 + \frac{e^2}{n - \beta}\right) - e^2 - 2 \beta (n - \beta) \Delta^2\right).$$
If 
$$Inv(G) = \frac{n - \beta}{\Delta} + \frac{1}{\beta \Delta^3} \left( 2 \beta \left(\beta \delta^2 + \frac{e^2}{n - \beta}\right) - e^2 - 2 \beta (n - \beta) \Delta^2\right).$$
then, from the above proofs, we have that  $\sum_{u \in I} d(u) = e$ and $\sum_{v \in V - I} d(v) = e$ and thereby $G$ is a
bipartite graph with partition sets of $I$ and $V - I$. Furthermore, we have that $d(u) = \delta = \Delta$ for each $u \in I$ and  
$d(v) = \Delta$ for each $v \in I$. Thus $\Delta |I| = |E(I, V - I)| = (n - |I|) \Delta$. Therefore $|I| = |V - I| = \frac{n}{2}$.
Hence $G$ is a regular balanced bipartite graph. \\

\noindent If $G$ is a regular balanced bipartite graph, a simple computation can verify that 
$$Inv(G) = \frac{n - \beta}{\Delta} + \frac{1}{\beta \Delta^3} \left( 2 \beta \left(\beta \delta^2 + \frac{e^2}{n - \beta}\right) - e^2 - 2 \beta (n - \beta) \Delta^2\right).$$

This completes the proofs of [$5$]. \\

\noindent {\bf Proof of Theorem $2$.} Let $G$ be a $k$-connected ($k \geq 2$) graph with $n \geq 3$ vertices and $e$ edges. Suppose $G$ is not Hamiltonian.
Then Lemma $2$ implies that $\beta \geq k + 1$. Also, we have that $n \geq 2 \delta + 1 \geq 2 k + 1$ otherwise 
$\delta \geq k \geq n/2$
and $G$ is Hamiltonian.
 Let $I_1 := \{\, u_1, u_2, ..., u_{\beta} \,\}$ be a maximum independent set in $G$. Then
  $I := \{\, u_1, u_2, ..., u_{k + 1} \,\}$ is an independent set in $G$.
  Thus
$$ \sum_{u \in I} d(u) = |E(I, V - I)| \leq \sum_{v \in V - I} d(v).$$ 
Since $\sum_{u \in I} d(u) +  \sum_{v \in V - I} d(v) = 2e$, we have that 
$$\sum_{u \in I} d(u) \leq e \leq \sum_{v \in V - I} d(v).$$
\noindent [$1$]. Following the proof of [$1$] in Theorem $1$, we have 
$$Z_1(G) \leq  (n - k - 1)\Delta^2 + \frac{e^2}{2(k + 1)} + \frac{(k + 1) \Delta^3}{2 \delta}.$$
Since $$Z_1(G) \geq  (n - k - 1)\Delta^2 + \frac{e^2}{2(k + 1)} + \frac{(k + 1) \Delta^3}{2 \delta},$$
we have $$Z_1(G) =  (n - k - 1)\Delta^2 + \frac{e^2}{2(k + 1)} + \frac{(k + 1) \Delta^3}{2 \delta}.$$
Thus we further have that  $\sum_{u \in I} d(u) = e$ which implies that $\sum_{v \in V - I} d(v) = e$ and thereby $G$ is a
bipartite graph with partition sets of $I$ and $V - I$. Furthermore, we have that $d(u) = \delta = \Delta$ for each $u \in I$ and  
$d(v) = \Delta$ for each $v \in I$. Thus $\delta |I| = |E(I, V - I)| = (n - |I|) \delta$. Therefore $(k + 1) = |I| = |V - I| = \frac{n}{2}$.
By Lemma $4$, we have $G$ is Hamiltonian, a contradiction. \\

\noindent [$2$]. Following the proofs of [$2$] in Theorem $1$ and [$1$] in Theorem $2$, we can show that $G$ is Hamiltonian, a contradiction. 
The details of the proofs are skipped here. \\
 
\noindent [$3$]. Following the proofs of [$3$] in Theorem $1$ and [$1$] in Theorem $2$, we can show that $G$ is Hamiltonian, a contradiction. 
The details of the proofs are skipped here. \\

\noindent [$4$]. Following the proofs of [$4$] in Theorem $1$ and [$1$] in Theorem $2$, we can show that $G$ is Hamiltonian, a contradiction. 
The details of the proofs are skipped here. \\
 
\noindent [$5$]. Following the proofs of [$5$] in Theorem $1$ and [$1$] in Theorem $2$, we can show that $G$ is Hamiltonian, a contradiction. 
The details of the proofs are skipped here. \\

This completes the proofs of Theorem $2$. \\

\noindent {\bf Proof of Theorem $3$.} Let $G$ be a $k$-connected ($k \geq 1$) graph with $n \geq 9$ vertices and $e$ edges. Suppose $G$ is not traceable.
Then Lemma $3$ implies that $\beta \geq k + 2$. Also, we have that $n \geq 2 \delta + 2 \geq 2 k + 2$ otherwise 
$\delta \geq k \geq (n - 1)/2$
and $G$ is traceable.
 Let $I_1 := \{\, u_1, u_2, ..., u_{\beta} \,\}$ be a maximum independent set in $G$. Then
  $I := \{\, u_1, u_2, ..., u_{k + 2} \,\}$ is an independent set in $G$.
  Thus
$$ \sum_{u \in I} d(u) = |E(I, V - I)| \leq \sum_{v \in V - I} d(v).$$ 
Since $\sum_{u \in I} d(u) +  \sum_{v \in V - I} d(v) = 2e$, we have that 
$$\sum_{u \in I} d(u) \leq e \leq \sum_{v \in V - I} d(v).$$
\noindent [$1$]. Following the proof of [$1$] in Theorem $1$, we have 
$$Z_1(G) \leq  (n - k - 2)\Delta^2 + \frac{e^2}{2(k + 2)} + \frac{(k + 2) \Delta^3}{2 \delta}.$$
Since $$Z_1(G) \geq  (n - k - 2)\Delta^2 + \frac{e^2}{2(k + 2)} + \frac{(k + 2) \Delta^3}{2 \delta},$$
we have $$Z_1(G) =  (n - k - 2)\Delta^2 + \frac{e^2}{2(k + 2)} + \frac{(k + 2) \Delta^3}{2 \delta}.$$
Thus we further have that  $\sum_{u \in I} d(u) = e$ which implies that $\sum_{v \in V - I} d(v) = e$ and thereby $G$ is a
bipartite graph with partition sets of $I$ and $V - I$. Furthermore, we have that $d(u) = \delta = \Delta$ for each $u \in I$ and  
$d(v) = \delta$ for each $v \in I$. Thus $\delta |I| = |E(I, V - I)| = (n - |I|) \delta$. Therefore $(k + 2) = |I| = |V - I| = \frac{n}{2}$.
Since $n \geq 9$, we have $k \geq 3$. By Lemma $4$, we have $G$ is Hamiltonian and thereby $G$ is traceable, a contradiction. \\

\noindent [$2$]. Following the proofs of [$2$] in Theorem $1$ and [$1$] in Theorem $3$, we can show that $G$ is traceable, a contradiction. 
The details of the proofs are skipped here. \\
 
\noindent [$3$]. Following the proofs of [$3$] in Theorem $1$ and [$1$] in Theorem $3$, we can show that $G$ is traceable, a contradiction. 
The details of the proofs are skipped here. \\

\noindent [$4$]. Following the proofs of [$4$] in Theorem $1$ and [$1$] in Theorem $3$, we can show that $G$ is traceable, a contradiction. 
The details of the proofs are skipped here. \\
 
\noindent [$5$]. Following the proofs of [$5$] in Theorem $1$ and [$1$] in Theorem $3$, we can show that $G$ is traceable, a contradiction. 
The details of the proofs are skipped here. \\
 
 This completes the proofs of Theorem $3$. \\

\end{document}